# NONPARAMETRIC ESTIMATION OF A POINT-SPREAD FUNCTION IN MULTIVARIATE PROBLEMS

BY PETER HALL[1] AND PEIHUA QIU[2]

*Australian National University and University of Minnesota*

The removal of blur from a signal, in the presence of noise, is readily accomplished if the blur can be described in precise mathematical terms. However, there is growing interest in problems where the extent of blur is known only approximately, for example in terms of a blur function which depends on unknown parameters that must be computed from data. More challenging still is the case where no parametric assumptions are made about the blur function. There has been a limited amount of work in this setting, but it invariably relies on iterative methods, sometimes under assumptions that are mathematically convenient but physically unrealistic (e.g., that the operator defined by the blur function has an integrable inverse). In this paper we suggest a direct, noniterative approach to nonparametric, blind restoration of a signal. Our method is based on a new, ridge-based method for deconvolution, and requires only mild restrictions on the blur function. We show that the convergence rate of the method is close to optimal, from some viewpoints, and demonstrate its practical performance by applying it to real images.

**1. Introduction.** Observed signals are usually not exactly the same as true signals, but are instead degraded. This can occur through the entire process of signal acquisition, for a variety of reasons. For example, in aerial reconnaissance, astronomy and remote sensing, signals are often adversely affected by atmospheric turbulence or aberrations of the optical system. Signal degradations can be classified into several categories, among which *point degradation* (or noise) and *spatial degradation* (or blur) are the most common. Other types of degradation involve chromatic or temporal effects.

Received December 2005; revised September 2006.
[1]Supported in part by an ARC grant.
[2]Supported in part by an NSF grant and an NSA grant.

For a detailed account of the formation and nature of degradations the reader is referred to books such as those by Andrews and Hunt [1] and Bates and McDonnell [2]. Related discussion is also given by Qiu [18].

In image analysis the true signal is often observed, or scanned, on a two-dimensional pixel grid, subject to both noise and blur. More generally, a signal may be recorded in any number of dimensions. For example, Lidar imaging devices record in $d = 3$ dimensions, and a great deal of signal analysis is conducted in the case $d = 1$.

In these settings it can be considered that we observe

$$(1.1) \qquad Y(j) = (\phi\psi)(j) + N(j),$$

where $\psi$ denotes the true signal, $N$ represents noise, $\phi$ is a linear operator applied to $\psi$ and $Y$ is the noisy signal. The latter is acquired on a $d$-variate square lattice, and therefore $j$, in (1.1), is a member of the set $\mathbb{Z}^d$ of all $d$-vectors of integers. We shall use the symbol $\phi$ to denote also the kernel of the operator $\phi$; this function is sometimes referred to as the blur function. Thus, $(\phi\psi)(j) = \sum_k \phi(j-k)\,\psi(k)$, for each $j \in \mathbb{Z}^d$.

In an image-analysis interpretation of (1.1), $\psi$ denotes the true scene, $Y$ is the observed image, the function $\phi$ is called a point-spread function, $\mathbb{Z}^d$ is a mathematical representation of the pixel grid on the Charge Coupled Device (CCD) and of course, $d = 2$. In this setting, and also more generally, we expect $\phi$ to preserve signal intensity, that is, $\sum_j \phi(j) = 1$. In particular, this implies that if $\psi \equiv b$ for a constant $b$, meaning that the true signal is of constant "brightness," then $\phi\psi \equiv \psi$.

Image restoration (when $d = 2$), or, more generally, signal restoration, is a process for reconstructing a close approximation to the unobserved signal $\psi$ from its observed but degraded form, $Y$. Many procedures for image restoration assume that $\phi$ is known. This is the case with, for example, the inverse filter, Wiener filter, constrained least-squares filter, Lucy–Richardson procedure, Landweber procedure, Tikhonov–Miller procedure, maximum a posteriori (MAP) procedure, maximum entropy procedure and techniques based on the EM algorithm. See, for instance, [20], [11], Chapter 5, [4] and [10]. In some settings this is reasonable, since $\phi$ can be specified, at least approximately, using our knowledge of the signal acquisition device. However, in other applications this information is not available, and so approximation (or estimation) of $\phi$ is a prerequisite for image restoration. This is the context of the present paper.

Signal restoration when $\phi$ is unknown is referred to as blind signal restoration. A number of procedures have been proposed for solving this problem. They can be grouped into two categories. In the first, $\phi$ is described by a parametric model, usually with just one, but occasionally two, parameters. See, for example, the work of Cannon [3], Katsaggelos and Lay [15], Rajagopalan and Chaudhuri [19], Carasso [5] and Joshi and Chaudhuri [14].

The other class of procedures assumes that the true signal consists of an object with specific, known features—for example, a shape with known support, against a uniform background—but involves only weak, nonparametric assumptions about $\phi$. See, for example, the contributions of Yang, Galatsanos and Stark [22] and Kundur and Hatzinakos [16].

It is to the latter category that we contribute in this paper. We introduce a method which, working from a known test signal and making only mild, nonparametric assumptions about the blur function, recovers the latter without suffering the drawbacks of earlier nonparametric techniques. In particular, the mechanism leading to the observed signal is not precisely known because we lack information about the blur function, rather than about the true signal.

Using our technique, the blur function does not need to have an integrable inverse, or reciprocal. The latter assumption will very seldom be satisfied in practice, although it is made in recent, related literature. Moreover, our technique is substantially less complex than the iterative approaches which are invariably used in nonparametric settings.

We introduce a new, ridge-based deconvolution algorithm. Unlike conventional methods, this technique is well suited to inversion when the Fourier transformation of the point-spread function vanishes at infinitely many points. Standard approaches to dealing with this problem sometimes resort to "fencing off" those zeros, and then dealing separately with each one. That can be particularly awkward, and is avoided by our ridge-based method. In addition to having good numerical performance, the ridge technique achieves, in some settings, theoretically optimal convergence rates, and so is no less "sharp" than its more conventional competitors.

Our theoretical work is related to earlier contributions of Hall [12], Johnstone and Silverman [13], Donoho and Low [7], Donoho [6], Van Rooij, Ruymgaart and van Zwet [21] and Ermakov [8]. These authors, in a variety of settings, discuss consistency and convergence rates for inverse estimators computed from unknown signals and from known blur functions.

The method is introduced in Section 2. Some of its theoretical properties are discussed in Section 3. A numerical study in Section 4 describes our method's statistical features and its application to real images. Technical details are deferred to Section 5.

## 2. Models and estimators.

2.1. *Model for degraded, noisy signal.* We assume model (1.1) throughout. The noise, $N$, is taken to be independent and identically distributed at each lattice point $j$, with variance $\sigma^2 > 0$. We suppose that $\phi$ preserves signal intensity.

We might think of the pixel-based blur function $\phi$ as representing a discrete approximation to an idealised, smooth blur function, $g$ say, which operates in the continuum. If the pixel width is considered to be $n^{-1}$, where $n \geq 1$ is an integer which we shall permit to become arbitrarily large, then the relationship between $\phi$ and $g$ might be taken to be

$$
\begin{aligned}
\phi(j) &= n^{-d} f(j/n) \qquad \text{for all } j \in \mathbb{Z}^d \quad \text{and} \\
f(x) &= s_d g(x) \qquad \text{for all } x \in \mathbb{R}^d,
\end{aligned}
\tag{2.1}
$$

where we might take the function $g$ to be fixed (i.e. not depending on $n$) and, if $\phi$ preserves intensity (e.g., preserves the light energy striking the CCD in a typical imaging device), the scaling factor $s_d$ satisfies

$$
s_d = \left\{ n^{-d} \sum_{j \in \mathbb{Z}^d} g(j/n) \right\}^{-1} \to 1
\tag{2.2}
$$

as $n \to \infty$. The limiting relation in (2.2) holds because $\int g = 1$, this being the continuum version of $\sum_j \phi(j) = 1$. Thus, $f$ is a normalized version of $g$, on the pixel grid, and $\phi$ is a discretized version of $f$.

The suggestion that $g$ be a fixed function is made here only to simplify our ideas. In our subsequent theoretical work we shall, through analogous changes to $\phi$, permit the spread of $g$ to alter with $n$, so that the difficulty of the imaging problem can evolve as the amount of information changes.

We shall take $f$ to be supported on the sphere of radius $\lambda_n/n$. It follows that $g$ is supported on the same set.

2.2. *Model for test signal.* In the case $d = 2$, test signals, or test patterns, are frequently used to determine a point-spread function from data. Test patterns are images that are known to significantly greater accuracy than that provided by the image recording device under test. In fact, test patterns are generally known completely; there is no need to estimate parameters, and in this sense the term "parametric image model" would be misleading if it were applied to a test pattern in a narrow statistical sense. In practice, performance is often assessed visually; in this paper we use mathematical closeness in the $L_2$ metric in lieu of subjective assessment.

Real test patterns are typically comprised of regular geometric shapes, such as rectangles. We shall treat such a signal here, in the $d$-variate case, although to simplify notation and discussion we shall assume that there is a single rectangular prism, $m_j$ pixels wide along the $j$th axis for $j = 1, \ldots, d$. If the sides of the prism are parallel to the pixel axes, if the lower left- and upper right-hand corners of the rectangle are at $(a_1, \ldots, a_d)$ and $(b_1, \ldots, b_d)$, respectively, and if the value of the signal is 1 within the rectangular prism and 0 outside, then $\psi(k_1, \ldots, k_d)$ equals 1 if $a_j \leq k_j \leq b_j$ for $1 \leq j \leq d$, and equals zero otherwise. It follows that, with $t = (t_1, \ldots, t_d)^{\mathrm{T}}$ and $j = (j_1, \ldots, j_d)^{\mathrm{T}}$, we have

$$
\begin{aligned}
\psi^{\mathrm{Ft}}(t) &= \sum_{j_1=a_1}^{b_1} \cdots \sum_{j_d=a_d}^{b_d} e^{it^{\mathrm{T}}j} \\
&= \exp\left\{\frac{1}{2}i\sum_{\ell=1}^{d}(a_\ell + b_\ell)t_\ell\right\} \prod_{\ell=1}^{d} \frac{\sin(m_\ell t_\ell/2)}{\sin(t_\ell/2)},
\end{aligned}
\tag{2.3}
$$

where $m_\ell = b_\ell - a_\ell + 1$, and the superscript Ft denotes the discrete Fourier transform.

If, as in the discussion of (2.1) and (2.2) in Section 2.1, we consider the lattice $\mathbb{Z}^d$ to represent a rescaled pixel grid where neighbors are, in reality, distant $n^{-1}$ rather than 1 apart along each axis, then it is reasonable to consider $m_\ell$ to be asymptotic to $c_\ell n$, where in this instance we take $c_\ell > 0$ to be fixed as $n$ diverges. In this way the rectangular $m_1 \times \cdots \times m_d$ prism represents, as $n$ diverges and scale is suitably adjusted, an increasingly accurate approximation to a prism with edge lengths $c_1, \ldots, c_d$.

2.3. *Discrete Fourier transforms.* Assume that $\phi$ vanishes outside a known sphere $\mathcal{R} = \mathcal{R}(n)$ in $\mathbb{Z}^d$, centred at the origin, $O$, and of radius $\lambda_n$, where $n/\lambda_n$ is bounded; and that $\psi$ likewise is zero outside a known set $\mathcal{S}$, which extends no further than radius $O(n)$ from $O$. Put $\mathcal{T} = \mathcal{R} \oplus \mathcal{S} = \{j + k : \in \mathcal{R}, k \in \mathcal{S}\}$. Then $\phi\psi$ vanishes outside $\mathcal{T}$, and

$$
\phi^{\mathrm{Ft}}(t) = \sum_{j \in \mathbb{Z}^d} \phi(j)e^{it^{\mathrm{T}}j} = \sum_{j \in \mathcal{R}} \phi(j)e^{it^{\mathrm{T}}j}, \qquad \psi^{\mathrm{Ft}}(t) = \sum_{j \in \mathcal{S}} \psi(j)e^{it^{\mathrm{T}}j}
$$

and $(\phi\psi)^{\mathrm{Ft}} = \phi^{\mathrm{Ft}}\psi^{\mathrm{Ft}}$.

In a slight abuse of notation we denote by $Y^{\mathrm{Ft}}$ and $N^{\mathrm{Ft}}$ the Fourier transforms of $Y$ and $N$ restricted to $\mathcal{T}$,

$$
Y^{\mathrm{Ft}}(t) = \sum_{j \in \mathcal{T}} Y(j)e^{it^{\mathrm{T}}j}, \qquad N^{\mathrm{Ft}}(t) = \sum_{j \in \mathcal{T}} N(j)e^{it^{\mathrm{T}}j}.
$$

Therefore, a Fourier-transform version of (1.1) has the form

$$
Y^{\mathrm{Ft}}(t) = \phi^{\mathrm{Ft}}(t)\psi^{\mathrm{Ft}}(t) + N^{\mathrm{Ft}}(t), \qquad t \in \mathbb{R}^d.
\tag{2.4}
$$

Result (2.4) highlights the symmetry of the problem: In principle, identical methods can be used to recover $\phi$ from $Y$ knowing $\psi$, and to recover $\psi$ knowing $\phi$. However, a marked degree of asymmetry is often introduced through the typical forms of $\phi$ and $\psi$. Again the problem of image analysis provides a convenient example. There, when the point-spread function $\phi$ is known, and the problem is that of estimating the true scene, then $\phi$ is generally smooth, and in particular $\phi^{\mathrm{Ft}}(t)$ generally converges relatively quickly to zero as $\|t\|$ increases. (Here, $\|\cdot\|$ denotes the Euclidean metric on $\mathbb{R}^d$.) On the other hand, when the true scene is known, for example a test pattern, and the problem is one of estimating the point-spread function, $\psi$ is often unsmooth. In particular, as indicated in Section 2.2, $\psi$ contains jump discontinuities, representing the sharp boundaries in a test pattern. In such cases, $\psi^{\mathrm{Ft}}(t)$ generally converges to zero relatively slowly. Of course, there are exceptions to these generalities; for example, if $\phi$ denotes the point-spread function corresponding to motion blur then it is unsmooth.

We shall concentrate on the problem of estimating $\phi$ from known $\psi$.

2.4. *Estimation of $\phi$ from known $\psi$.* Let $\rho(t)$ denote a positive constant multiple of a known, positive function of the real variable, $t$. We use $\rho(t)$ as, in effect, a ridge when regularizing a Fourier transform. In particular, recognizing that $\phi^{\mathrm{Ft}} = (\phi\psi)^{\mathrm{Ft}}/\psi^{\mathrm{Ft}}$ and therefore

$$\phi(j) = \frac{1}{(2\pi)^d} \int_{\mathcal{A}} \frac{(\phi\psi)^{\mathrm{Ft}}(t)}{\psi^{\mathrm{Ft}}(t)} e^{-ij^{\mathrm{T}}t} \, dt, \tag{2.5}$$

where $\mathcal{A} = [-\pi, \pi]^d$, we define an estimator $\hat{\phi}$ of $\phi$ by

$$\hat{\phi}(j) = \frac{1}{(2\pi)^d} \int_{\mathcal{A}} \frac{\psi^{\mathrm{Ft}}(-t)|\psi^{\mathrm{Ft}}(t)|^r Y^{\mathrm{Ft}}(t)}{\{|\psi^{\mathrm{Ft}}(t)| \vee \rho(t)\}^{r+2}} e^{-ij^{\mathrm{T}}t} \, dt. \tag{2.6}$$

Here, $r \geq 0$; choosing $r > 0$ removes potential numerical problems associated with computing the integral in (2.6).

We may think of (2.6) as having been obtained from (2.5) by (a) multiplying the numerator and the denominator in the integrand of (2.5) by $\psi^{\mathrm{Ft}}(-t)|\psi^{\mathrm{Ft}}(t)|^r$; (b) replacing $|\psi^{\mathrm{Ft}}(t)|$ by the maximum of that quantity and the ridge, in the quantity $|\psi^{\mathrm{Ft}}(t)|^{r+2}$ which step (a) produces in the denominator; and (c) replacing $(\phi\psi)^{\mathrm{Ft}}(t)$ in the numerator by its unbiased approximation, $Y^{\mathrm{Ft}}(t)$.

In some cases, considerations of symmetry in the process for manufacturing the signal recording device imply that, to a first approximation, $\phi$ is radially symmetric. For example, glass (as distinct from resin) lens elements are typically manufactured using a polishing process which involves rolling a large sphere, with cylindrical glass blanks attached, inside another sphere. However, errors in this process can introduce asymmetric aberrations to such elements, in particular because the outer, grinding sphere is worn, or the glass blanks are not correctly secured. Other causes of asymmetry result from inaccuracies in the alignment of elements within the lens, or in the cementing of lens elements together. Since the design of a lens is often highly complex, there are many different ways in which asymmetric aberrations can arise, and no standard parametric models for the blur functions that they might produce.

## 3. Theoretical properties.

3.1. *Mean-squared error criteria for choosing the ridge.* We define the sum of squared errors of $\hat{\phi}(j)$ to be

$$\text{SSE} = \sum_{j \in \mathbb{Z}^d} |\hat{\phi}(j) - \phi(j)|^2.$$

From this formula and the definition of $\hat{\phi}(j)$, at (2.6), it follows that the mean summed squared error (MSSE) admits the formula

$$\text{MSSE} = E(\text{SSE}) = \sigma^2(\#\mathcal{T})\frac{1}{(2\pi)^d}\int_{\mathcal{A}}\left\{\frac{|\psi^{\text{Ft}}|^{r+1}}{(|\psi^{\text{Ft}}| \vee \rho)^{r+2}}\right\}^2$$

(3.1)
$$+ \frac{1}{(2\pi)^d}\int_{\mathcal{A}}|\phi^{\text{Ft}}|^2\left\{1 - \frac{|\psi^{\text{Ft}}|^{r+2}}{(|\psi^{\text{Ft}}| \vee \rho)^{r+2}}\right\}^2.$$

The first and second terms on the right-hand side represent the total contributions to mean summed squared error from variance and squared bias, respectively. For example, the first term on the far right-hand side equals

$$\sum_{j \in \mathbb{Z}^d} E|\hat{\phi}(j) - E\hat{\phi}(j)|^2 = \frac{1}{(2\pi)^d}\int_{\mathcal{A}}\left\{\frac{|\psi^{\text{Ft}}|^{r+1}}{(|\psi^{\text{Ft}}| \vee \rho)^{r+2}}\right\}^2 E|N^{\text{Ft}}|^2,$$

and the claimed relationship follows from the fact that $E|N^{\text{Ft}}|^2 \equiv \sigma^2(\#\mathcal{T})$.

We suggested in Section 2.1 that the lattice $\mathbb{Z}^d$ be interpreted as a rescaled version of a pixel grid with edge length $n^{-1}$. We claim that in this setting it is appropriate to work with $n^d\text{MSSE}$, rather than directly with MSSE. To appreciate why, recall from (2.1) and (2.2) that $n^d\phi(j)$ can be interpreted as a discrete approximation, on a pixel grid, to a continuous blur function $f$ evaluated at $j/n$. In this context, $\hat{f}(j/n) \equiv n^d\hat{\phi}(j)$ can be viewed as an estimator of $f(j/n)$ and extended to $\mathbb{R}^d$; and $n^d\text{MSSE}$ can be interpreted as a discrete approximation to the mean integrated squared error of $\hat{f}$ as an approximation to $f$.

3.2. *Asymptotic properties of $\phi^{\text{Ft}}$ and $\psi^{\text{Ft}}$.* Reflecting the rescaling discussed above, define $\phi_n^{\text{Ft}}(t) = \phi^{\text{Ft}}(t/n)$, $\psi_n^{\text{Ft}}(t) = n^{-d}\psi^{\text{Ft}}(t/n)$, $\rho_n = n^{-d}\rho$, $\mathcal{A}_n = [-n\pi, n\pi]^d$ and $\tau = n^{-d}(\#\mathcal{T})$. In this notation,

$$n^d\text{MSSE} = n^{-d}\sigma^2\tau\frac{1}{(2\pi)^d}\int_{\mathcal{A}_n}\left\{\frac{|\psi_n^{\text{Ft}}|^{r+1}}{(|\psi_n^{\text{Ft}}| \vee \rho_n)^{r+2}}\right\}^2$$

(3.2)
$$+ \frac{1}{(2\pi)^d}\int_{\mathcal{A}_n}|\phi_n^{\text{Ft}}|^2\left\{1 - \frac{|\psi_n^{\text{Ft}}|^{r+2}}{(|\psi_n^{\text{Ft}}| \vee \rho_n)^{r+2}}\right\}^2.$$

Since $n^d\text{MSSE}$ can be represented so simply in terms of $\phi_n^{\text{Ft}}$ and $\psi_n^{\text{Ft}}$, then it is of interest to know properties of those functions.

We shall work with classes of compactly-supported blur functions $\phi$ for which the associated, rescaled Fourier transform, $\phi_n^{\text{Ft}}$, decreases at least polynomially fast in the tails. In general it is awkward to prove that such a rate of decrease occurs arbitrarily far out in the tails, but fortunately we need it only a distance $o(n)$ from the origin.

Performance is determined by three main parameters: $n$, the number of observations per linear unit of space; $\sigma_n^2$, noise variance; and $\lambda_n$, blur radius. Arguably the first two of these are the most intrinsic, although $\lambda_n$ also plays a major role.

These considerations lead us to define the following class of sequences of blur functions. Given a sequence $\Lambda = \{\lambda_n\}$ of positive numbers, and $p > 0$:

Let $\mathcal{F}(\Lambda, p)$ denote a class of sequences of functions $\phi$ depending on $n$, with the following properties for each given $n$: (a) $\phi$ vanishes outside a $d$-variate sphere, centered at the origin, of radius $\lambda_n$; (b) $\sum_j \phi(j) = 1$;

(3.3) and (c) for each positive sequence $\varepsilon_n$ decreasing to zero as $n \to \infty$,

$(1 + \|t\|)^p |\phi_n^{\text{Ft}}(t)|$ is bounded uniformly in $t \in \mathcal{A}_n$ with $\|t\| \leq n\varepsilon_n$

and in $\phi \in \mathcal{F}(\Lambda, p)$, and $|\phi_n^{\text{Ft}}(t)| = O\{(n\varepsilon_n)^{-d}\}$ uniformly in $t \in \mathcal{A}_n$

with $\|t\| > n\varepsilon_n$ and in $\phi \in \mathcal{F}(\Lambda, p)$.

Here and below, "$\phi \in \mathcal{F}(\Lambda, p)$" means that the sequence of blur functions for which the function, at "time" $n$, is $\phi$, is in $\mathcal{F}(\Lambda, p)$. Thus, $\phi$ depends on the pixel scale-factor $n$, although to prevent ambiguity we indicate this in notation only for the Fourier transform $\phi_n^{\text{Ft}}$ of the rescaled version of $\phi$, not for $\phi$ itself.

To illustrate, we introduce a function $\phi$ which satisfies the conditions in (3.3). The function class $\mathcal{F}(\Lambda, p)$ could be taken to be a set of rescaled versions of this $\phi$, but of course it can be much larger.

Consider a compactly supported, continuum blur function, $g(x) = A_1 \prod_\ell (1 - x_\ell^2)^p$ for $\sup_\ell |x_\ell| \leq 1$. The constant $A_1 > 0$ is chosen so that the function integrates to 1 on $[-1, 1]^d$, or equivalently, so that signal intensity is preserved. The associated characteristic function,

$$g^{\text{Ft}}(t) = A_1 \int_{x : \sup_\ell |x_\ell| \leq 1} e^{it^{\text{T}} x} \left\{ \prod_{\ell=1}^{d} (1 - x_\ell^2)^p \right\} dx,$$

satisfies $|g^{\text{Ft}}(t)| \leq A_2 (1 + \|t\|)^{-p}$, for all $t \in \mathbb{R}^d$, where $A_2 > 0$ is a constant.

The discrete blur function $\phi$ analogous to $g$ is

(3.4) $$\phi(j) = A_2(n) n^{-d} \prod_{\ell=1}^{d} (1 - n^{-2} j_\ell^2)^p$$

for $\sup_\ell |j_\ell| \leq n$, where the bounded sequence $A_2(n)$ is chosen to preserve signal intensity. In this case, $\lambda_n = O(n)$; that is, $\phi$ is supported within a sphere of radius $n$ of the origin. For each sequence $\varepsilon_n \downarrow 0$ there exists a constant $A_3 > 0$ such that the corresponding $\phi_n^{\mathrm{Ft}}$ satisfies $|\phi_n^{\mathrm{Ft}}(t)| \leq A_3(1+\|t\|)^{-p}$ for all $\|t\| \leq n\varepsilon_n$, and also $|\phi_n^{\mathrm{Ft}}(t)| \leq A_3(n\varepsilon_n)^{-d}$.

The signal model introduced at (2.3) admits a concise asymptotic description. Let us, in (2.3), take $c_\ell = m_\ell/n$; then

$$(3.5) \qquad \psi_n^{\mathrm{Ft}}(t) = e^{\theta_n^{\mathrm{T}} t} \prod_{\ell=1}^d \frac{\sin(c_\ell t_\ell/2)}{n\sin(t_\ell/2n)},$$

where $\theta_n \in \mathbb{R}^d$. If $c_\ell$ is either fixed or converges to a finite, nonzero number as $n \to \infty$, then $\psi_n^{\mathrm{Ft}}(t)$ is asymptotic to $\psi_{\lim}^{\mathrm{Ft}}(t) = \prod_\ell \{2t_\ell^{-1}\sin(c_\ell t_\ell/2)\}$. It follows that, for each sequence $\varepsilon_n \downarrow 0$, $|\psi_n^{\mathrm{Ft}}(t) - \psi_{\lim}^{\mathrm{Ft}}(t)|/|\psi_{\lim}^{\mathrm{Ft}}(t)| \to 0$ uniformly in $\|t\| \leq n\varepsilon_n$, and $|\psi_n^{\mathrm{Ft}}(t)| = O\{(n\varepsilon_n)^{-d}\}$ uniformly in $t \in \mathcal{A}_n$ for which $\|t\| > n\varepsilon_n$.

3.3. *Upper bound to rate of convergence of MSSE.* Our main result is the following. Define $\mathcal{F}(\Lambda, p)$ as at (3.3). The formula for the threshold, $\rho_n(t)$, given there can be changed without appreciably altering the results. For example, the theorem continues to hold if, in the expression for $\rho_n(t)$, we replace $\prod_\ell(|t_\ell| \vee 1)$ by simply 1, strengthen the condition $p > \frac{1}{2}(d+q)$ to $p > d + \frac{1}{2}q$ and weaken the assumption $q > 3d$ to $q > 2d$.

THEOREM 3.1. *Assume that $n/\lambda_n$ is bounded as $n \to \infty$, and that the noise variance, $\sigma^2 = \sigma_n^2$, satisfies $n^{-d}\lambda_n^d \sigma_n^2 \to 0$ as $n \to \infty$. Let $h_n$ denote a positive sequence decreasing to zero, put $\rho_n(t) = h_n\{\prod_\ell(|t_\ell| \vee 1)\}^{-1}\|t\|^q$, take $r \geq 0$ in the definition of $\hat{\phi}(j)$ at (2.6) and assume that $p > \frac{1}{2}(d+q)$ and $q > 3d$. Then, as $n \to \infty$,*

$$(3.6) \qquad \sup_{\phi \in \mathcal{F}(\Lambda,p)} n^d \mathrm{MSSE} = O\{(n^{-d}\lambda_n^d \sigma_n^2 h_n^{-1} + h_n)(\log n)^{d-1}\}.$$

REMARK 3.1 (*Optimizing choice of $h_n$*). The theorem implies that a mean-square convergence rate of essentially $(\lambda_n^d \sigma_n^2/n^d)^{1/2}$, uniformly over $\phi \in \mathcal{F}(\Lambda, p)$, can be achieved by taking $h_n \asymp (\lambda_n^d \sigma_n^2/n^d)^{1/2}$:

$$(3.7) \qquad \sup_{\phi \in \mathcal{F}(\Lambda,p)} n^d \mathrm{MSSE} = O\{(\lambda_n^d \sigma_n^2/n^d)^{1/2}(\log n)^{d-1}\}.$$

The notation $a_n \asymp b_n$, for positive $a_n$ and $b_n$, means that $a_n/b_n$ is bounded away from zero and infinity as $n \to \infty$.

REMARK 3.2 (*Smoothness of $\phi$*). The convergence rate in (3.7) does not depend on the smoothness of $\phi$, represented by $p$ in the function class $\mathcal{F}(\Lambda, p)$, provided $p$ exceeds $\frac{1}{2}(d+q)$. It is of interest to consider what this means in terms of the number of derivatives enjoyed by the blur functions. Let us take $q = 3d+$, that is, just a little larger than $3d$. Then the condition $p > \frac{1}{2}(d+q)$ reduces to $p > 2d$. If $\mathcal{F}(\Lambda, p)$ is sufficiently large, for example if it contains a scale-changed version of the $\phi$ defined at (3.4), then the assumption that all the blur functions in $\mathcal{F}(\Lambda, p)$ have $s$ square-integrable derivatives is equivalent to asking that $p > s + \frac{1}{2}d$. In this setting the smoothness condition imposed in the theorem reduces to the restriction that all the functions in the class have $d$ bounded derivatives. In the important special case of image analysis, $d = 2$ and just two derivatives are required.

REMARK 3.3 (*Smoothness of $\psi$*). If the test signal, $\psi$, is a relatively smooth function, then the mean-square accuracy of even an optimal approximation to $\phi$ can be inferior to the rate in (3.7). For example, taking $\lambda_n = n$, $\sigma_n^2 = n^{-1}$ and $d = 1$ for simplicity, the rate in (3.7) is $n^{-1/2}$. However, assuming that $|\psi_n^{\mathrm{Ft}}|$ decreases like $(1 + \|t\|)^{-s}$ as $\|t\|$ diverges, the minimax-optimal rate of convergence of mean-squared error in estimation of $\phi$ can be shown to equal $n^{-2p/(2p+2s+1)}$. (See [9] for related results in density deconvolution problems.) This is inferior to the rate $n^{-1/2}$ unless $p > s + \frac{1}{2}$. Therefore, if the test signal is very smooth, the blur function must be even smoother if the accuracy of the estimator of the blur function is not to be degraded relative to that for a rough test signal.

3.4. *Lower bound to rate of convergence of MSSE.* Let $f$ denote a fixed, spherically symmetric, compactly supported probability density on $\mathbb{R}^d$, for which

$$\sup_{t \in \mathbb{R}^d} (1 + \|t\|)^p |f^{\mathrm{Ft}}(t)| < \infty, \tag{3.8}$$

where $p > 0$. Let $\xi$ be a $d$-vector, and put $\delta_n = \lambda_n^{-1}$,

$$\chi_\theta(x) = c_{1,\theta} \delta_n^d f(\delta_n x) \{1 + \theta \cos(\xi^{\mathrm{T}} x)\}, \tag{3.9}$$

where $\theta = 0$ or $1$ and $c_{1,\theta}$ denotes a constant. (Here, and below, we suppress the dependence of quantities such as $\chi_\theta$ and $c_{1,\theta}$ on $n$.) Note that

$$\delta_n^d \int f(\delta_n x) \cos(\xi^{\mathrm{T}} x) e^{it^{\mathrm{T}} x} \, dx$$
$$= \frac{1}{2} \left\{ f^{\mathrm{Ft}}\left(\frac{\xi + t}{\delta_n}\right) + f^{\mathrm{Ft}}\left(\frac{\xi - t}{\delta_n}\right) \right\}.$$

Therefore, if we define $c_{1,\theta}^{-1} = 1 + \theta f^{\mathrm{Ft}}(\xi/\delta_n)$, then $\chi_\theta$ is a probability density.

Let the blur function $\phi_\theta$ denote the conventional discrete approximation to $\chi_\theta$,

$$\phi_\theta(j) = c_{2,\theta} n^{-d} \chi_\theta(j/n), \qquad j \in \mathbb{Z}, \tag{3.10}$$

where the constant $c_{2,\theta}$ is chosen so that intensity is preserved, that is, $\sum_j \phi_\theta(j) = 1$. Under the conditions given in the theorem below, this standardization entails $c_{2,\theta} \to 1$ as $n \to \infty$. Let $\mathrm{MSSE}_{n\theta}$ denote the version of MSSE for a general estimator of $\phi$ [not necessarily the estimator at (2.6)], when the true $\phi$ is $\phi_\theta$ and the scale parameter equals $n$.

THEOREM 3.2. *Assume $\delta_n \asymp (\sigma_n^2/n^d)^{1/(3d)}$, that $C_1 n^{-C_2} \leq \delta_n^{2d} \leq C_3 n^{-C_4}$, where $C_1, \ldots, C_4 > 0$ and $d(1 - p^{-1}) \leq C_4 \leq C_2$, and that $p \geq \max(\frac{3}{2}d, \frac{3}{4}C_2 + 1)$. Suppose, too, that the noise variables $N(j)$ are independent and identically distributed as Normal $\mathrm{N}(0, \sigma_n^2)$. Then, for a choice of $\xi$ in (3.9) that depends only on $c_1, \ldots, c_d$ in the definition of the test signal $\psi$ [see (3.5)],*

$$(3.11) \qquad \liminf_{n\to\infty}(n^d/\lambda_n^d \sigma_n^2)^{1/2} \sup_{\theta=0,1} n^d \mathrm{MSSE}_{n\theta} > 0.$$

In view of (3.8), the sequence of functions $\phi_\theta$, indexed by $n$, is $\mathcal{F}(\Lambda, p)$ for $\theta = 0, 1$, provided the constant in the uniform bound on $(1 + \|t\|^p)|\phi_n^{\mathrm{Ft}}(t)|$ in (3.3) is chosen sufficiently large. In this case, (3.11) implies that

$$(3.12) \qquad \liminf_{n\to\infty}(n^d/\sigma_n^2)^{1/2} \sup_{\phi \in \mathcal{F}(\{\delta_n\}, p)} n^d \mathrm{MSSE}_{n\theta} > 0.$$

Assuming the relation $\lambda_n \asymp (n^d/\sigma_n^2)^{1/(3d)}$ between noise variance and support of the blur function, and with the exception of the logarithmic factor in (3.7), (3.12) is a converse of (3.7). Within these constraints, the estimator $\hat{\phi}$ at (2.6) recovers $\phi$ from the test-pattern data at the optimal rate.

In the case $\sigma_n^2 = n^{-1}$, treated in Remark 3.3, Theorem 3.2 shows that the convergence rate achieved is optimal if $\lambda_n \asymp n^{(d+1)/(3d)}$. This is a more realistic assumption than $\lambda_n \asymp n$, imposed in Remark 3.3, since it permits the number of pixels that represent the width of the blur function to be an order of magnitude less than the number per linear unit of space.

## 4. Numerical results.

4.1. *Square-block test pattern.* Here we summarize the results of a simulation study when $d = 2$, in cases where the true image, represented by the function $\psi$, is a simple square block with intensity 1, against a white background with intensity 0. See panel (a) of Figure 1. We take the true continuum blurring function to be

$$(4.1) \qquad g(x_1, x_2) = \frac{1}{(0.7388\lambda)^2}\{1 - (x_1/\lambda)^2\}\{1 - (x_2/\lambda)^2\}$$

for $\sup(|x_1|, |x_2|) \leq \lambda$, and $g(x_1, x_2) = 0$ otherwise, as suggested in Section 3.2 with $p = 5$. Note that $g$ is not circularly symmetric; that is, $g(x_1, x_2)$ is not a function of $x_1^2 + x_2^2$ alone. (The great majority of parametric models for point-spread functions are circularly symmetric.) We denote the discretized version of $g$ by $\phi$.

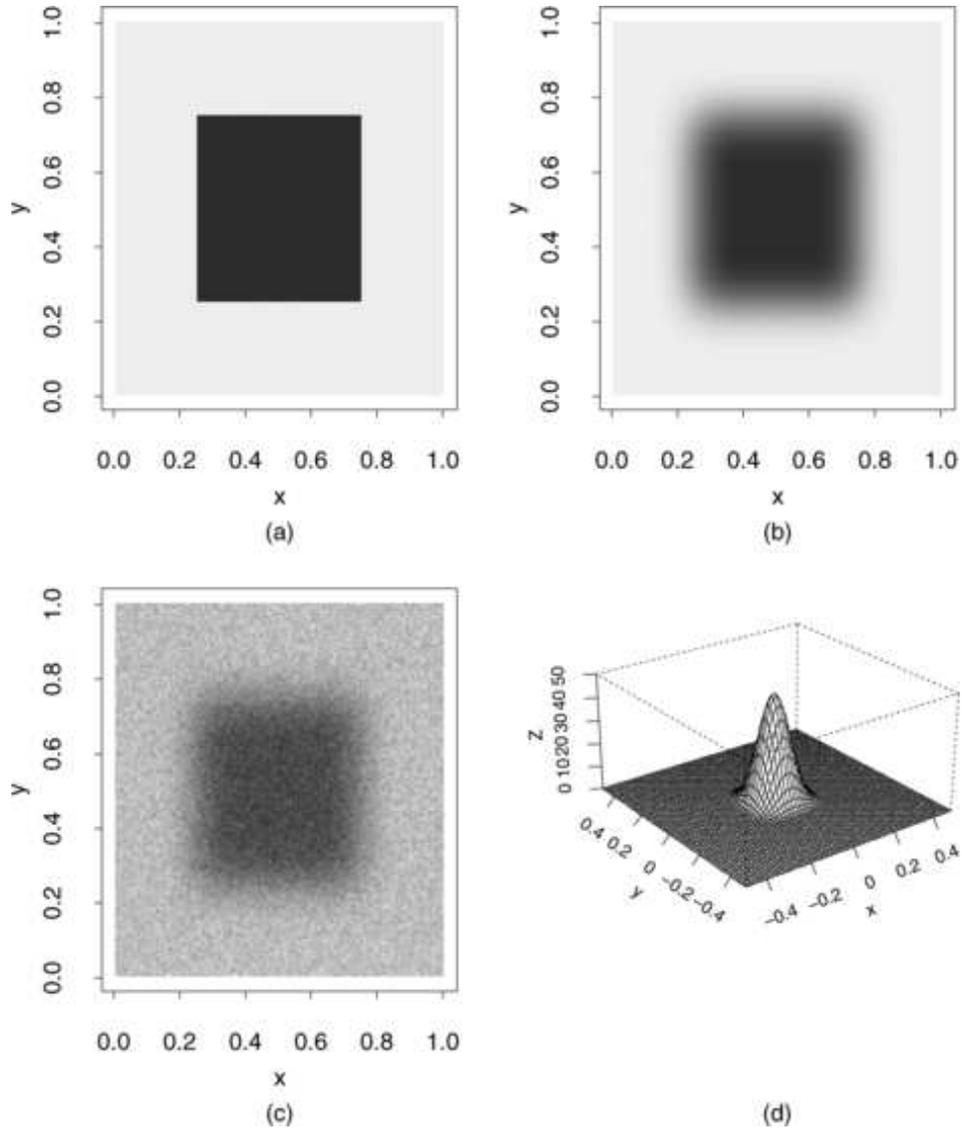

FIG. 1. *Graphs of $\psi$, $\phi$, $\phi\psi$ and $Y$. Panel* (a) *shows the true image, $\psi$. Panel* (b) *shows its blurred version, $\phi\psi$. The function $\phi$ itself is depicted in panel* (d). *Panel* (c) *shows the blurred image plus noise; the latter was $\mathrm{N}(0, 0.1^2)$ on each pixel. Digitization was on a $128 \times 128$ grid.*

TABLE 1
*Mean summed squared error. Tabulated are values of MSSE (the first number in each entry) and se of SSE (the second number) of the estimator $\hat{\phi}$ defined at (2.6), based on 101 simulations. The optimal value of $h_n$ is presented in the second line in each entry. The noise distribution is $N(0, \sigma^2)$*

| $r$ | $\sigma = 0.05$ | $\sigma = 0.1$ | $\sigma = 0.2$ |
|---|---|---|---|
| 1 | 1.9744, 0.0607 | 2.0987, 0.0750 | 2.3457, 0.0887 |
| | $1.0 \times 10^{-4}$ | $1.1 \times 10^{-4}$ | $1.2 \times 10^{-4}$ |
| 10 | 1.2046, 0.0779 | 1.3708, 0.0474 | 1.4961, 0.0954 |
| | $2.0 \times 10^{-5}$ | $3.1 \times 10^{-5}$ | $3.2 \times 10^{-5}$ |
| 50 | 0.6397, 0.0324 | 0.6644, 0.0531 | 0.7533, 0.0874 |
| | $1.7 \times 10^{-5}$ | $1.7 \times 10^{-5}$ | $1.7 \times 10^{-5}$ |
| 55 | 0.6397, 0.0324 | 0.6643, 0.0531 | 0.7532, 0.0873 |
| | $1.7 \times 10^{-5}$ | $1.7 \times 10^{-5}$ | $1.7 \times 10^{-5}$ |

See Figure 1(d) for a perspective plot of $\phi$ when $n = 128$ and $\lambda = 0.2$. Figure 1(b) shows the result of blurring $\psi$ using $\phi$. If we add independent and identically distributed $N(0, \sigma^2)$ noise to the blurred image at each pixel, then we obtain, when $\sigma = 0.1$, the result shown in Figure 1(c).

We evaluated the performance of the estimator $\hat{\phi}$, defined at (2.6), when $n = 128$, $\lambda = 0.2$ and $\sigma = 0.05$, 0.1 and 0.2. For the estimator $\hat{\phi}$ we chose $\rho_n(t) = h_n \|t\|^5$, which, along with $g$ in (4.1), satisfies the conditions given in Section 3.2. There are two parameters, $r$ and $h_n$, involved in the estimator $\hat{\phi}$. We found that, in most cases (e.g., $n = 128$ or 256, $\sigma \in [0, 1]$), results were improved when $r$ increased in the range $[0, 50]$, and they did not change much when $r$ was chosen larger than 50. However, when $r$ was chosen too large (e.g., larger than 60), numerical underflow sometimes occurred in the computations, since the denominator in (2.6) was very small in such cases. To demonstrate this, we consider four $r$ values: 1, 10, 50 and 55. For each combination of $\sigma$ and $r$, we searched for the optimal value of $h_n$ in the range $[0, 10^{-3}]$, with step-length $10^{-5}$. In this analysis we employed MSSE (mean summed squared error) to define optimality, as in Section 3, and used as our data the results of 101 simulations. Values of MSSE, the standard error of SSE and the optimal value of $h_n$ are presented in Table 1.

From Table 1 it can be seen that: (a) In all cases considered, MSSE values are stable when $r$ is chosen larger than 50; (b) MSSE increases with $\sigma$, but the effect of $\sigma$ is quite small; (c) $h_n$ should be chosen smaller when $r$ is larger or $\sigma$ is smaller.

We found that, for the smaller sample sizes treated in our numerical work, the estimator performed well except that it under-estimated the peak of $\phi$ a little. This is a common aberration of nonparametric curve and surface estimators, which tend to be biased down in peaks and up on troughs. The tendency can be largely removed by making a simple change of scale,

(4.2) $$\bar{\phi}(x_1, x_2) = \hat{\phi}(x_1/s, x_2/s),$$

where $s > 0$ is a tuning parameter.

In practice, all tuning parameters, including $s$, would be chosen to give the best visual impression. This approach is common in image analysis, and avoids difficulties that arise when using mathematical criteria that are based on $L_2$ performance but do not approximate visual perception particularly well. See [17] for discussion.

Figure 2(a) shows the estimator $\bar{\phi}$ that has median value of MSSE, out of 101 simulations, when $\sigma = 0.1$, $r = 50$, $h_n = 1.7 \times 10^{-5}$ and $s = 0.92$. Its profiles in the cross-sections of $x_2 = 0$ and $x_1 = 0$ are shown in Figures 2(b) and 2(c), respectively, by the dotted curves. In these two plots, the solid curves denote the profiles of the true point-spread function $\phi$, and the short and long-dashed curves denote the profiles of the estimator $\hat{\phi}$ having median value of MSSE, out of 101 simulations, when $r = 50$ and $\sigma = 0.05$, 0.1 and 0.2, respectively.

4.2. *Application to cameraman image.* To illustrate how our methodology affects the restored image in the entire image restoration process, we used the popular cameraman image as an example. The original image is shown in Figure 3(a); it is of size $256 \times 256$ pixels, with gray levels in the range $[0, 255]$. A blurred version of this image, using the point-spread function $g$ at (4.1) with $\lambda = 0.05$ (i.e. with a $25 \times 25$ pixel blurring window) is shown in Figure 3(b). Figure 3(c) depicts the image that is obtained after adding independent and identically distributed $N(0, 5^2)$ noise to the image in Figure 3(b).

We pretended that these images were made by the same image acquisition device as that for the test image shown in Figure 1. Then the point-spread function, $\phi$, was estimated by (2.6) and (4.2) from the degraded test image, using the same level of blurring and noise as for the degraded cameraman images. We fixed $r$ at 50, as before.

There are several existing procedures for restoring $\psi$ from $Y$, if $\phi$ is known or estimated. We chose two noniterative procedures: the inverse filter with a hard threshold, and the Wiener filter. The restored image computed by the first approach is given by

$$\hat{\psi}_1(x) = \frac{1}{(2\pi)^2} \Re\bigg\{\int\!\!\int \frac{Y^{\mathrm{Ft}}(t)}{\hat{\phi}^{\mathrm{Ft}}(t)} I(|\hat{\phi}^{\mathrm{Ft}}(t)| > \gamma) \exp(it^{\mathrm{T}} x)\, dt\bigg\},$$

where $\hat{\phi}$ denotes the estimated point-spread function, and $\gamma > 0$ is the threshold. The restored image obtained by the second approach is defined by

$$\hat{\psi}_2(x) = \frac{1}{(2\pi)^2} \Re\bigg\{\int\!\!\int \frac{\hat{\phi}_*^{\mathrm{Ft}}(t)}{|\hat{\phi}^{\mathrm{Ft}}(t)|^2 + \alpha\|t\|^\beta} Y^{\mathrm{Ft}}(t) \exp(it^{\mathrm{T}} x)\, dt\bigg\},$$

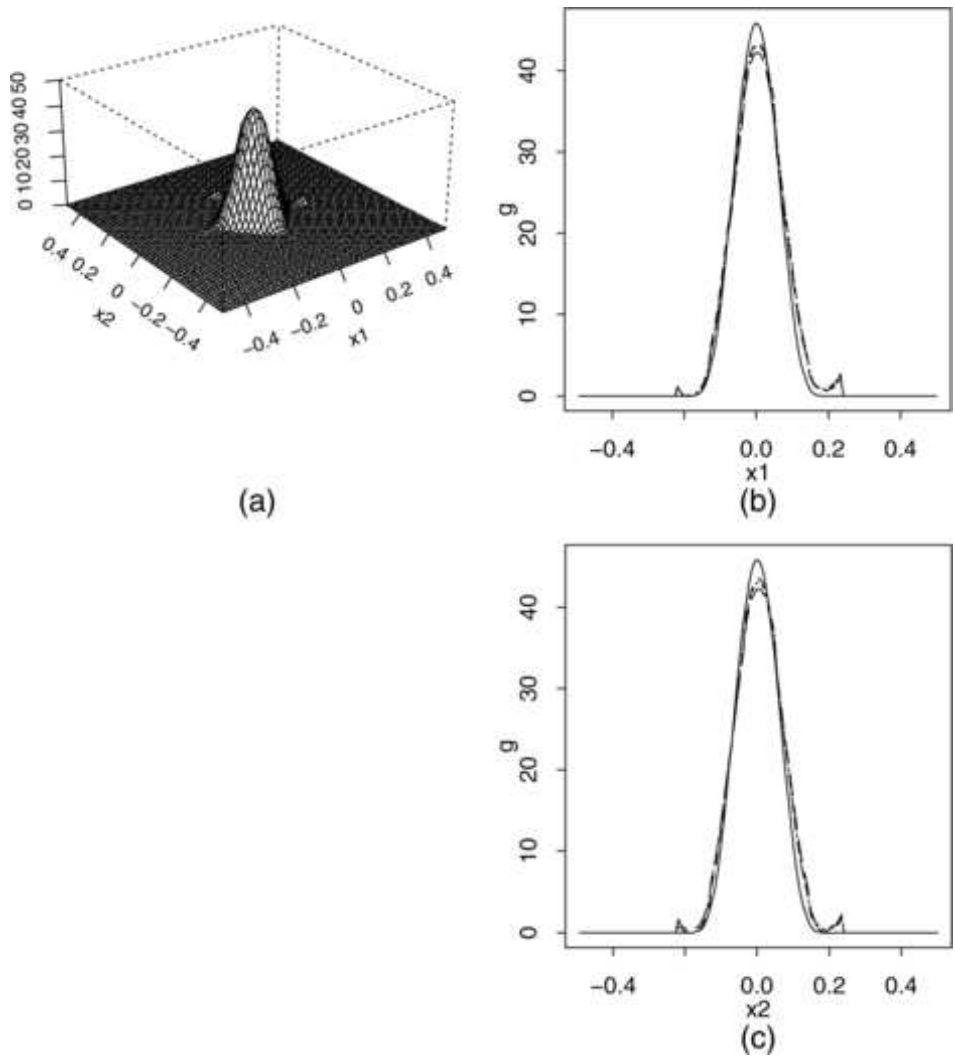

FIG. 2. *Graphs of $\bar{\phi}$. Panel* (a) *shows a plot of $\bar{\phi}$ when $\sigma = 0.1$, $r = 50$, $h_n = 1.7 \times 10^{-5}$ and $s = 0.92$. Panels* (b) *and* (c) *show profiles in the cross-sections of $x_2 = 0$ and $x_1 = 0$, respectively, of $\phi$ (solid), $\bar{\phi}$ when $\sigma = 0.05$ (dotted), $\bar{\phi}$ when $\sigma = 0.1$ (short-dashed) and $\bar{\phi}$ when $\sigma = 0.2$ (long-dashed). In each case, the estimator $\bar{\phi}$ has median value of MSSE, out of 101 simulations.*

where $\hat{\phi}_*^{\mathrm{Ft}}$ denotes the complex conjugate of $\hat{\phi}^{\mathrm{Ft}}$, and $\alpha, \beta > 0$ are two parameters. The inverse filter is basically the least-squares procedure; use of the threshold alleviates noise amplification. The Wiener filter is derived with a view to minimizing MSSE of the restored image under the assumption that noise is Gaussian. These two approaches are used commonly in the literature. See [11], Chapter 5, for detailed discussion.

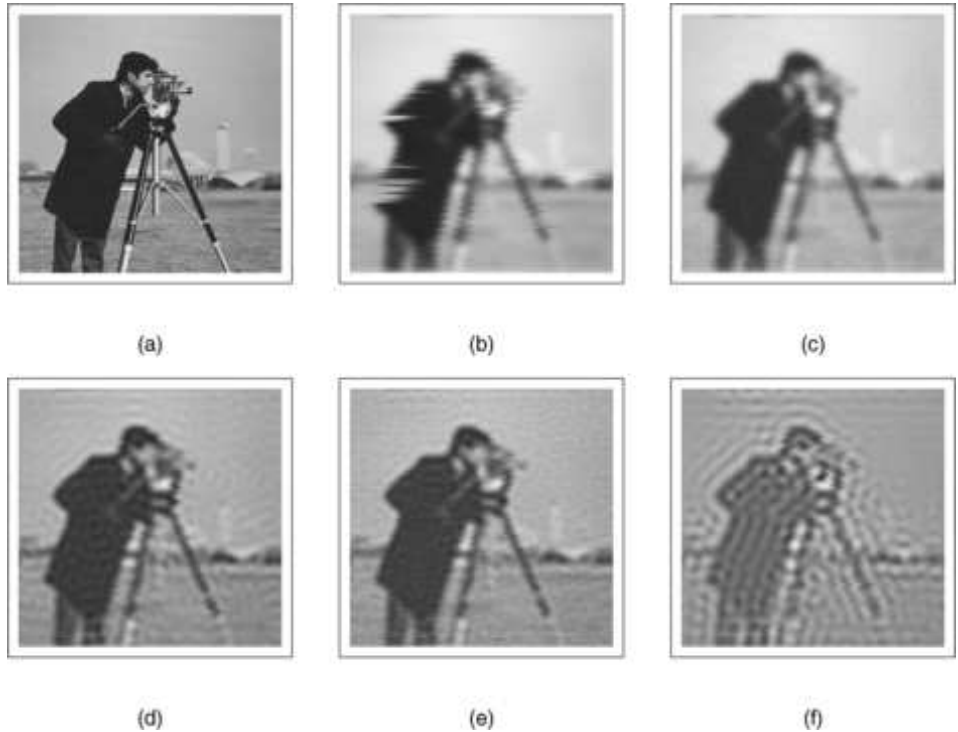

FIG. 3. *Cameraman example. Panels* (a)–(c) *show the original, the blurred, and the blurred-and-noisy cameraman images, respectively. Panels* (d) *and* (e) *show images restored from* (c) *by the inverse filter and the Wiener filter, respectively. Panel* (f) *shows the restored image, obtained by the Wiener filter, when a Gaussian point-spread function with standard deviation* $\lambda/2$ *was used in deblurring.*

The restored image, obtained by inverse filtering from the blurred and noisy cameraman image, is shown in Figure 3(d). The corresponding results for Wiener filtering are given in Figure 3(e). In each case the tuning parameters, $\alpha$, $\beta$, $h_n = 10^{-5}$ and $s = 0.89$, were selected to give a good visual impression.

Next, instead of estimating the point-spread function as suggested at (2.6), we assumed that the Wiener filter used a Gaussian point-spread function with its standard deviation equal to $\lambda/2$. (This produces virtually the best results in the Gaussian case. Note that the radius of the Gaussian point-spread function is effectively twice its standard deviation, and that the radius of the correct point-spread function equals $\lambda$.) The corresponding result is shown in Figure 3(f). It can be seen that this mistaken guess at the point-spread function affects the results considerably.

## 5. Proofs.

5.1. *Proof of Theorem* 3.1. Define $\beta_1(t) = \prod_\ell |t_\ell^{-1} \sin(c_\ell t_\ell/2)|$, with $c_\ell$ as in (3.5); put

$$u_{n1} = \int_{\mathcal{A}_n} \left[ \frac{|\psi_n^{\mathrm{Ft}}(t)|^{r+1}}{\{|\psi_n^{\mathrm{Ft}}(t)| \vee \rho_n(t)\}^{r+2}} \right]^2 dt,$$

$$v_{n1} = \int_{\mathcal{A}_n} |\phi_n^{\mathrm{Ft}}(t)|^2 \left[ 1 - \frac{|\psi_n^{\mathrm{Ft}}(t)|^{r+2}}{\{|\psi_n^{\mathrm{Ft}}(t)| \vee \rho_n(t)\}^{r+2}} \right]^2 dt;$$

let $\gamma$ denote a general, positive function of $t = (t_1, \ldots, t_d)$ that depends on the $t_j$'s only through their absolute values; and let $u_{n2}(\gamma)$ and $v_{n2}(\gamma)$ have the same respective definitions as $u_{n1}$ and $v_{n1}$, but with $|\psi_n^{\mathrm{Ft}}(t)|$ and $\rho_n(t)$ replaced by $\beta_1(t)$ and $\rho_n(t)\gamma(t)$, respectively. Noting that the denominator contribution to $\psi_n^{\mathrm{Ft}}(t)$ at (3.5) satisfies

$$A \prod_{\ell=1}^d |t_\ell| \leq \prod_{\ell=1}^d |n \sin(t_\ell/2n)| \leq \tfrac{1}{2} \prod_{\ell=1}^d |t_\ell|,$$

uniformly in $t \in \mathcal{A}_n$, where $A \in (0, \frac{1}{2})$ is an absolute constant, and with $w$ denoting either $u$ or $v$, we have

(5.1) $$w_{n1} \leq \sup_\gamma{}' w_{n2}(\gamma),$$

uniformly in $\phi \in \mathcal{F}(\Lambda, p)$, where $\sup'_\gamma$ denotes the supremum over choices of $\gamma$ satisfying $B_1^{-1} \leq \gamma \leq B_1$.

Let $\alpha(t) = (1 + \|t\|)^{-p}$, $\beta_2(t) = \prod_\ell |t_\ell^{-1} \sin(t_\ell/2)|$,

(5.2) $$u_{n3}(\gamma) = \int_{\mathcal{A}_n} \left[ \frac{\beta_2(t)^{r+1}}{\{\beta_2(t) \vee \rho_n(t)\gamma(t)\}^{r+2}} \right]^2 dt,$$

(5.3) $$v_{n3}(\varepsilon, \gamma) = \int_{\|t\| \leq n\varepsilon} |\alpha(t)|^2 \left[ 1 - \frac{\beta_2(t)^{r+2}}{\{\beta_2(t) \vee \rho_n(t)\gamma(t)\}^{r+2}} \right]^2 dt,$$

where $0 < \varepsilon < \pi$. If we change variables in the integrals defining $u_{n2}$ and $v_{n2}$, from $t = (t_1, \ldots, t_d)^{\mathrm{T}}$ to $s = (s_1, \ldots, s_d)^{\mathrm{T}}$, with $s_\ell = c_\ell t_\ell$ where $c_\ell$ is as in (3.5), and if we observe that, in the definition of $u_{n2}(\gamma)$, the method for bounding the integral over a rectangle $\prod_\ell [nc_{1\ell}, nc_{2\ell}]$, where $-\infty < c_{1\ell} < 0 < c_{2\ell} < \infty$, is the same as that for the integral over $\mathcal{A}_n$, then it can be deduced from (3.2), (5.1), the fact that $\#\mathcal{T} = O(\lambda_n^d)$, and the definition of $\mathcal{F}(\Lambda, p)$ that, for each positive sequence $\varepsilon_n$ decreasing to zero, and for $B_1, B_2 > 0$ sufficiently large,

(5.4) $$n^d \mathrm{MSSE} \leq B_2 \sup_\gamma{}' \{n^{-d} \lambda_n^d \sigma_n^2 u_{n3}(\gamma) + v_{n3}(\varepsilon_n, \gamma)\} + O(n^{-d} \lambda_n^d \sigma_n^2 \varepsilon_n^{-2d}),$$

uniformly in $\phi \in \mathcal{F}(\Lambda, p)$.

Define
$$I_d(\varepsilon) = \int_{\mathcal{D}_d} \beta(s)^z I\{\beta(s) \le \varepsilon\}\, ds,$$
where $\varepsilon > 0$, $z \ge 0$, $\mathcal{D}_d = [0,1]^d$ and $\beta(s) = \prod_{1 \le \ell \le d} |s_\ell|$. The result below describes the size of $I_d(\varepsilon)$. □

LEMMA 5.1. *As $\varepsilon \downarrow 0$, $I_d(\varepsilon) = O(\varepsilon^{z+1} |\log \varepsilon|^{d-1})$.*

PROOF. Observe that
$$I_d(\varepsilon) = \int_{\mathcal{D}_{d-1}} \left(\prod_{\ell=2}^d s_\ell\right)^z I\left(\prod_{\ell=2}^d s_\ell \le \varepsilon\right) ds_2 \cdots ds_d \int_0^1 s_1^z\, ds_1$$
(5.5)
$$+ \int_{\mathcal{D}_{d-1}} \left(\prod_{\ell=2}^d s_\ell\right)^z I\left(\prod_{\ell=2}^d s_\ell > \varepsilon\right) ds_2 \cdots ds_d \int_0^{\varepsilon/s_2 \cdots s_d} s_1^z\, ds_1$$
$$= I_{d-1}(\varepsilon) + \frac{\varepsilon^{z+1}}{z+1} J_{d-1}(\varepsilon),$$
where
$$J_d(\varepsilon) = \int_{\mathcal{D}_d} I\left(\prod_{\ell=1}^d s_\ell > \varepsilon\right) \frac{ds_1 \cdots ds_d}{s_1 \cdots s_d}$$
$$= \int_{\mathcal{D}_{d-1}} \log\left(\frac{s_2 \cdots s_d}{\varepsilon}\right) I\left(\prod_{\ell=2}^d s_\ell > \varepsilon\right) \frac{ds_2 \cdots ds_d}{s_2 \cdots s_d}$$
$$\le |\log \varepsilon| J_{d-1}(\varepsilon).$$

The latter inequality, and an argument by induction, imply that $J_d(\varepsilon) \le |\log \varepsilon|^d$. This bound and (5.5) establish that $I_d(\varepsilon) \le I_{d-1}(\varepsilon) + (z+1)^{-1} \varepsilon^{z+1} \times |\log \varepsilon|^{d-1}$. It is readily proved that $I_1(\varepsilon) = (z+1)^{-1} \varepsilon^{z+1}$, and so it follows inductively that $I_d(\varepsilon) = O(\varepsilon^{z+1} |\log \varepsilon|^{d-1})$ as $\varepsilon \downarrow 0$, completing the proof of the lemma. □

Next we give a bound for $v_{n3}(\varepsilon_n, \gamma)$, with $v_{n3}(\varepsilon, \gamma)$ defined as at (5.3). If $j = (j_1, \ldots, j_d)$, where each component is an integer, let the $d$-variate cube $\mathcal{C}_j$ denote the set of $t = (t_1, \ldots, t_d)$ for which each $t_\ell - j_\ell \pi \in [-\frac{1}{2}\pi, \frac{1}{2}\pi)$. Taking $\varepsilon = \varepsilon_n$ at (5.3), we may bound the integral there by the sum, $v_{n4}(\varepsilon_n, \gamma)$ say, over vectors $j$ for which $\|j\| \le 2n\varepsilon_n$, of the integrals
$$K_j \equiv \int_{\mathcal{C}_j} |\alpha(t)|^2 \left[1 - \frac{\beta_2(t)^{r+2}}{\{\beta_2(t) \vee \rho_n(t)\gamma(t)\}^{r+2}}\right]^2 dt.$$

In turn, $v_{n4}(\varepsilon_n, \gamma) = v_{n5}(\varepsilon_n, \gamma) + v_{n6}(\varepsilon_n, \gamma)$, where $v_{n5}(\varepsilon_n, \gamma)$ equals the sum of $K_j$ over the set $\mathcal{K}_n$ of indices $j$ for which each $|j_\ell| \geq 1$ and $\|j\| \leq 2n\varepsilon_n$. Below, we shall establish an order-of-magnitude bound for $\sup'_\gamma v_{n5}(\varepsilon_n, \gamma)$, uniformly in $\phi \in \mathcal{F}(\Lambda, p)$. Similar methods may be use to derive the same bound for $\sup'_\gamma v_{n6}(\varepsilon_n, \gamma)$.

Define $s_\ell = \frac{1}{2}(t_\ell - j_\ell \pi)$, $\mathcal{C} = [-\frac{1}{4}\pi, \frac{1}{4}\pi]^d$, $\mathcal{D} = \mathcal{D}_d = [0,1]^d$ and $\beta(s) = \prod_\ell |s_\ell|$, the latter for $s = (s_1, \ldots, s_d) \in \mathcal{C}$. Since $\rho_n(t) = h_n(\prod_\ell |t_\ell|)^{-1}\|t\|^q$, then, for each $j \in \mathcal{K}_n$ and $t \in \mathcal{C}_j$, $\beta_2(t) \vee \rho_n(t)\gamma(t) = \theta_j(t)\{\beta(s) \vee h_n \|j\|^q\}$, where $\theta_j(t)$ is bounded away from zero and infinity uniformly in such $j$ and $t$. Therefore, defining $\delta_n(u) = h_n u^q$ for $u > 0$, and assuming that $\varepsilon_n \to 0$ so slowly that $n\varepsilon_n \to \infty$, we have

$$\sup_\gamma{}' v_{n5}(\varepsilon_n, \gamma)$$

(5.6)
$$= O\left(\sum_{j \in \mathcal{K}_n}(1 + \|j\|)^{-2p}\int_\mathcal{C}\left[1 - \frac{\beta(s)^{r+2}}{\{\beta(s) \vee (h_n\|j\|^q)\}^{r+2}}\right]^2 ds\right)$$

$$= O\left(\int_1^{2n\varepsilon_n}(1+u)^{d-1-2p}\,du\int_\mathcal{D}\left[1 - \frac{\beta(s)^{r+2}}{\{\beta(s) \vee \delta_n(u)\}^{r+2}}\right]^2 ds\right),$$

uniformly in $\phi \in \mathcal{F}(\Lambda, p)$. With $\delta = \delta_n(u)$ we have, uniformly in $1 \leq u \leq 2n\varepsilon_n$,

(5.7)
$$\int_\mathcal{D}\left[1 - \frac{\beta(s)^{r+2}}{\{\beta(s) \vee \delta\}^{r+2}}\right]^2 ds = \int_\mathcal{D}[1 - \{\beta(s)/\delta\}^{r+2}]^2 I\{\beta(s) \leq \delta\}\,ds$$
$$\leq \int_\mathcal{D} I\{\beta(s) \leq \delta\}\,ds$$
$$\leq \text{const.}\delta(1 + |\log \delta|)^{d-1},$$

where the last inequality is a consequence of Lemma 5.1.

From (5.3) and (5.7) it follows that, provided $p > \frac{1}{2}(d + q)$,

$$\sup_\gamma{}' v_{n5}(\varepsilon_n, \gamma) = O\left\{h_n(\log n)^{d-1}\int_1^{2n\varepsilon_n}(1+u)^{d+q-1-2p}\,du\right\}$$
$$= O\{h_n(\log n)^{d-1}\},$$

uniformly in $\phi \in \mathcal{F}(\Lambda, p)$. An identical bound applies to $\sup'_\gamma v_{n6}(\varepsilon_n, \gamma)$, and therefore to $\sup'_\gamma v_{n4}(\varepsilon_n, \gamma)$ and so to $\sup'_\gamma v_{n3}(\varepsilon_n, \gamma)$,

(5.8)
$$\sup_\gamma{}' v_{n3}(\varepsilon_n, \gamma) = O\{h_n(\log n)^{d-1}\}.$$

A similar argument shows that, with $u_{n3}(\gamma)$ as at (5.2), $\mathcal{L}_n$ denoting the set of $j$ for which each $|j_\ell| \geq 1$ and $\|j\| \leq n\pi$, and $\langle j \rangle = (\prod_{1 \leq \ell \leq d}|j_\ell|)^2$, we

have

$$\sup_{\gamma}{}' u_{n3}(\gamma) = O\left[\sum_{j \in \mathcal{L}_n} \langle j \rangle \int_{\mathcal{D}} \left\{\frac{\beta(s)^{r+1}}{\delta_n(j)^{r+2}}\right\}^2 I\{\beta(s) \le \delta_n(j)\}\, ds\right]$$

$$= O\left[\sum_{j \in \mathcal{L}_n} \langle j \rangle \delta_n(j)^{-2(r+2)} \int_{\mathcal{D}} \beta(s)^{2(r+1)} I\{\beta(s) \le \delta_n(j)\}\, ds\right]$$

$$= O\left\{\sum_{j \in \mathcal{L}_n} \langle j \rangle \delta_n(j)^{-1} (\log n)^{d-1}\right\}$$

$$= O\left\{h_n^{-1} (\log n)^{d-1} \int_{\mathcal{B}_n} \left(\prod_{\ell=1}^d x_j\right)^2 \left(1 + \sum_{j=1}^d x_j^3\right)^{-q/3} dx\right\},$$

where $\mathcal{B}_n$ denotes the set of $x \in \mathbb{R}^d$ for which each $x_\ell \ge 0$ and $\|x\| \le n$, and the second-last relation follows from Lemma 5.1. Changing variable from $x_j$ to $y_j = x_j^3$ in the last-written integral, we see that the integral is uniformly bounded provided that $q > 3d$. In this case,

$$\sup_{\gamma}{}' u_{n3}(\gamma) = O\{h_n^{-1}(\log n)^{d-1}\}. \tag{5.9}$$

Combining (5.4), (5.8) and (5.9), we deduce that

$$n^d \mathrm{MSSE} = O\{(n^{-d}\lambda_n^d \sigma_n^2 h_n^{-1} + h_n)(\log n)^{d-1} + n^{-d}\lambda_n^d \sigma_n^2 \varepsilon_n^{-2d}\}. \tag{5.10}$$

Since $\varepsilon_n$ here can be taken to equal any sequence that converges to zero more slowly than $n^{-1}$, then the theorem follows from (5.10).

5.2. *Proof of Theorem* 3.2. Let $\sigma_n^2$ denote the noise variance, let $N_{(0,1)}$ be a random variable having the $N(0,1)$ distribution and define

$$r_n = \sigma_n^{-2} \sum_{j \in \mathbb{Z}^d} \{(\phi_0 \psi)(j) - (\phi_1 \psi)(j)\}^2.$$

Consider the problem of deciding between $\theta = 0$ and $\theta = 1$ on the basis of the data $Y(j)$, defined at (1.1), for $j \in \mathcal{T}$. This is a classification problem, for which the likelihood-ratio rule consists of deciding in favor of $\theta = 0$ if

$$\sum_{j \in \mathbb{Z}^d} \{Y(j) - (\phi_0 \psi)(j)\}^2 \le \sum_{j \in \mathbb{Z}^d} \{Y(j) - (\phi_1 \psi)(j)\}^2,$$

and deciding in favor of $\theta = 1$ otherwise. From this property it can be proved that

(5.11)     the probability that the likelihood-ratio rule decides for $\theta = 1$ when $\theta = 0$,

or for $\theta = 0$ when $\theta = 1$, equals $\pi_n \equiv P(2N_{(0,1)} > r_n^{1/2})$.

Define $\chi_{n\theta}^{\text{Ft}}(t) = n^{-d}\sum_j \chi_\theta(j/n)e^{it^{\text{T}}j/n}$. Using Parseval's identity and employing the argument leading to (3.2), it can be shown that

$$(2\pi)^d r_n = \frac{n^d}{\sigma_n^2} I_1, \tag{5.12}$$

where $I_1 = \int_{\mathcal{A}_n} |\phi_{n0}^{\text{Ft}} - \phi_{n1}^{\text{Ft}}|^2 |\psi_n^{\text{Ft}}|^2$, $\phi_{n\theta}^{\text{Ft}}(t) = \phi_\theta^{\text{Ft}}(t/n)$ and $\phi_\theta^{\text{Ft}}(t) = \sum_j \phi_\theta(j) e^{it^{\text{T}}j}$, with $\phi_\theta$ defined at (3.10). Using the Euler–Maclaurin summation formula it can be proved that

$$\sup_{t\in\mathcal{A}_n} |\chi_{n\theta}^{\text{Ft}}(t) - \chi_\theta^{\text{Ft}}(t)| = O(n^{1-p}). \tag{5.13}$$

Since $\chi_\theta^{\text{Ft}}(0) = 1$ and the definition of $c_{2,\theta}$ is equivalent to $c_{2,\theta}\chi_{n\theta}^{\text{Ft}}(0) = 1$, then (5.13) implies that $c_{2,\theta} = 1 + O(n^{1-p})$. Therefore, noting that $\phi_{n\theta}^{\text{Ft}}(t) = c_{2,\theta}\chi_{n\theta}^{\text{Ft}}(t)$, we see that (5.13) implies that

$$\sup_{t\in\mathcal{A}_n} |\phi_{n\theta}^{\text{Ft}}(t) - \chi_\theta^{\text{Ft}}(t)| = O(n^{1-p}).$$

This result, and the fact that $\int_{\mathcal{A}_n} |\psi_n^{\text{Ft}}|^2 = O(1)$, imply that

$$|I_1^{1/2} - I_2^{1/2}| = O(n^{1-p}), \tag{5.14}$$

where $I_2 = \int_{\mathcal{A}_n} |\chi_0^{\text{Ft}} - \chi_1^{\text{Ft}}|^2 |\psi_n^{\text{Ft}}|^2$.

Observe that $I_2 = \frac{1}{4}c_{1,1}^2 I_3$ and $|I_3^{1/2} - I_4^{1/2}| = O(b_1^2)$, where

$$I_3 = \int_{\mathcal{A}_n} \left| b_1 f^{\text{Ft}}(t/\delta_n) + f^{\text{Ft}}\left(\frac{\xi+t}{\delta_n}\right) + f^{\text{Ft}}\left(\frac{\xi-t}{\delta_n}\right) \right|^2 |\psi_n^{\text{Ft}}(t)|^2 \, dt,$$

$$I_4 = \int_{\mathcal{A}_n} \left| f^{\text{Ft}}\left(\frac{\xi+t}{\delta_n}\right) + f^{\text{Ft}}\left(\frac{\xi-t}{\delta_n}\right) \right|^2 |\psi_n^{\text{Ft}}(t)|^2 \, dt$$

and $b_1 = 2(1 - c_{1,1}^{-1})$. For the choice $\delta_n \asymp (\sigma_n^2/n^d)^{1/(3d)}$ that we shall ultimately make,

$$c_{1,1} = \{1 + f^{\text{Ft}}(\xi/\delta_n)\}^{-1} = 1 + O(\delta_n^p) = 1 + O(n^{-pC_4/(2d)}) = 1 + O(n^{(1-p)/2}),$$

where we have used the fact that $C_4 \geq d(1-p^{-1})$. Therefore, $|I_3^{1/2} - I_4^{1/2}| = O(n^{1-p})$, and so by (5.14),

$$|I_1^{1/2} - \tfrac{1}{2}c_{1,1} I_4^{1/2}| = O(n^{1-p}). \tag{5.15}$$

We shall assume that each $c_\ell$ in (3.5) equals 1; the contrary case can be treated by changing variable in each coordinate. Then

$$I_4 = \int_{\mathcal{A}_n} \left| f^{\mathrm{Ft}}\left(\frac{\xi+t}{\delta_n}\right) + f^{\mathrm{Ft}}\left(\frac{\xi-t}{\delta_n}\right) \right|^2 \left| \prod_{\ell=1}^{d} \frac{\sin(t_\ell/2)}{n \sin(t_\ell/(2n))} \right|^2 dt$$

(5.16) $$\leq \mathrm{const.} \int_{\mathcal{A}_n} \left\{ \left| f^{\mathrm{Ft}}\left(\frac{\xi+t}{\delta_n}\right) \right|^2 + \left| f^{\mathrm{Ft}}\left(\frac{\xi-t}{\delta_n}\right) \right|^2 \right\} \left| \prod_{\ell=1}^{d} \frac{\sin(t_\ell/2)}{t_\ell} \right|^2 dt$$

$$\equiv \mathrm{const.} I_5,$$

say. Take $\xi = (2\pi, \ldots, 2\pi)^{\mathrm{T}}$. Then $I_5$ can be decomposed into a sum of two integrals, of which the first is

$$I_6 = \int_{\mathcal{A}_n} \left| f^{\mathrm{Ft}}\left(\frac{\xi+t}{\delta_n}\right) \right|^2 \left| \prod_{\ell=1}^{d} \frac{\sin(t_\ell/2)}{t_\ell} \right|^2 dt$$

and the second we denote by $I_7$. We shall show how to bound $I_6$; $I_7$ can be treated similarly.

Let $\mathcal{A}_{n1}$ be the set of points in $t = (t_1, \ldots, t_d)^{\mathrm{T}} \in \mathcal{A}_n$ for which $|t_\ell - 2\pi| > \pi$ for some $\ell$, and put $\mathcal{A}_{n2} = \mathcal{A}_n \setminus \mathcal{A}_{n1}$. The contribution to $I_6$ from integrating over $\mathcal{A}_{n1}$ equals $O(\delta_n^{2p})$. To bound the contribution, say $I_8$, to $I_6$ from integrating over $\mathcal{A}_{n2}$, note that on the latter set, $\prod_\ell t_\ell$ is bounded above zero. Therefore, changing variable from $t$ to $s$ where $t = \delta_n s - \xi$, we obtain

$$I_8 \leq \mathrm{const.} \int_{\mathcal{A}_{n2}} \left| f^{\mathrm{Ft}}\left(\frac{\xi+t}{\delta_n}\right) \right|^2 \left| \prod_{\ell=1}^{d} \sin(t_\ell/2) \right|^2 dt$$

$$\leq \mathrm{const.} \delta_n^d \int_{\mathbb{R}^d} \left| f^{\mathrm{Ft}}(s) \prod_{\ell=1}^{d} \sin(\delta_n s_\ell/2) \right|^2 ds$$

$$\leq \mathrm{const.} \delta_n^{3d} \int_{\mathbb{R}^d} \left| f^{\mathrm{Ft}}(s) \prod_{\ell=1}^{d} s_\ell \right|^2 ds = O(\delta_n^{3d}),$$

the identity holding because $p > 3d/2$. Therefore, $I_6 = O(\delta_n^{3d})$, and an identical bound can be derived for $I_7$, implying that $I_5 = O(\delta_n^{3d})$, and hence, by (5.16), that $I_4 = O(\delta_n^{3d})$. Therefore, in view of (5.15), $I_1 = O(\delta_n^{3d} + n^{2-2p})$. Since $p \geq \frac{3}{4}C_2 + 1$, then, for the choice $\delta_n \asymp (\sigma_n^2/n^d)^{1/(3d)}$, and noting that $\delta_n^{2d} \geq C_1 n^{-C_2}$, we have $n^{2-2p} = O(\delta_n^{3p})$, and thus, $I_1 = O(\delta_n^{3d})$. Hence, by (5.12),

(5.17) $$r_n = O(n^d \sigma_n^{-2} \delta_n^{3d}).$$

Define

$$s_n^2 = n^d \sum_{j \in \mathbb{Z}^d} \{\phi_0(j) - \phi_1(j)\}^2,$$

$$I_9 = \int_{\mathcal{A}_n} \left| f^{\mathrm{Ft}}\left(\frac{\xi+t}{\delta_n}\right) + f^{\mathrm{Ft}}\left(\frac{\xi-t}{\delta_n}\right) \right|^2 dt.$$

Arguments similar to those leading to (5.15) imply that, for a constant $C > 0$,

$$|s_n - \{1 + o(1)\} C I_9^{1/2}| = O(n^{1-p}).$$

From this property and for the choice $\delta_n \asymp (\sigma_n^2/n^d)^{1/(3d)}$, noting that $I_9 \asymp \delta_n^d$, and also that $p > \frac{1}{4}C_2 + 1$ [which entails $n^{1-p} = o(\delta_n^d)$], it can be shown that

(5.18) $$s_n^2 \asymp \delta_n^d.$$

Write $P_\theta$ and $E_\theta$ for probability measure and expectation, respectively, when the true blur function is $\phi_\theta$. Let $\pi_n$ be as in (5.11), and let $\eta > 0$. Let $\eta_j > 0$ denote a positive quantity which depends on $\eta$ but always satisfies $0 < \eta_j < 1$. Result (5.17) implies that if

(5.19) $$n^d \sigma_n^{-2} \delta_n^{3d} \leq \eta,$$

then $\frac{1}{2} \leq \pi_n \leq \frac{1}{2}(1 + \eta_1)$. Hence, by (5.11) and the Neyman–Pearson lemma, if $\hat{\theta}_n$ is any data-determined rule for deciding between $\theta = 0$ and $\theta = 1$,

(5.20) $$\liminf_{n \to \infty} \{P_0(\hat{\theta}_n = 1) + P_1(\hat{\theta}_n = 0)\} \geq 1 - \eta_2.$$

For the given the estimator $\hat{\phi}$, define $\hat{\theta}_n = 0$ if

$$\sum_{j \in \mathbb{Z}^d} |\hat{\phi}(j) - \phi_0(j)|^2 \leq \sum_{j \in \mathbb{Z}^d} |\hat{\phi}(j) - \phi_1(j)|^2,$$

and put $\hat{\theta}_n = 1$ otherwise. Then

$$\mathrm{SSE}_{n\theta} = \sum_{j \in \mathbb{Z}^d} |\hat{\phi}(j) - \phi_\theta(j)|^2 \geq \tfrac{1}{4} I(\hat{\theta}_n \neq \theta) n^{-d} s_n^2,$$

where the inequality follows from the triangle inequality. Therefore,

$$\sup_{\theta=0,1} n^d \mathrm{MSSE}_{n\theta} = \sup_{\theta=0,1} n^d E_\theta(\mathrm{SSE}_{n\theta})$$

$$\geq \tfrac{1}{4} s_n^2 \sup_{\theta=0,1} P_\theta(\hat{\theta}_n \neq \theta)$$

$$\geq \tfrac{1}{8} s_n^2 \{P_0(\hat{\theta}_n = 1) + P_1(\hat{\theta}_n = 0)\}.$$

This result and (5.20) imply that there exists $B_1 > 0$ such that

(5.21) $$\liminf_{n \to \infty} s_n^{-2} \sup_{\theta=0,1} n^d \mathrm{MSSE}_{n\theta} \geq B_1.$$

If we choose $\delta_n \asymp (\sigma_n^2/n^d)^{1/(3d)}$, and such that $\delta_n \leq (\eta \sigma_n^2/n^d)^{1/(3d)}$, then (5.19) holds and, using (5.18) to get the first inequality, $s_n^{-2} \leq B_2 \delta_n^{-d} \leq B_3(\sigma_n^2/n^d)^{-1/3}$. It follows from this result and (5.21) that

$$\liminf_{n\to\infty} (n^d/\sigma_n^2)^{1/3} \sup_{\theta=0,1} n^d \mathrm{MSSE}_{n\theta} \geq B_1 B_3^{-1},$$

which implies (3.11).

Department of Mathematics and Statistics
The University of Melbourne
Melbourne, Victoria 3010
Australia
E-mail: halpstat@ms.unimelb.edu.au

University of Minnesota
313 Ford Hall
224 Church Street, SE
Minneapolis, Minnesota 55455
USA
E-mail: qiu@stat.umn.edu